\documentclass[11pt,a4paper]{article}
\usepackage[utf8]{inputenc}

\usepackage{amssymb,amscd,amsfonts,amsbsy,amsmath,amsthm,amsrefs}
\usepackage{dsfont}
\usepackage{enumerate}
\usepackage{mathrsfs}
\usepackage{epsf,epsfig,esint}
\usepackage{pdfsync}
\usepackage[all]{xy}
\usepackage{pdfpages}

\newtheorem{proposition}{Proposition}[section]
\newtheorem{theorem}[proposition]{Theorem}
\newtheorem{lemma}[proposition]{Lemma}

\theoremstyle{definition}
\newtheorem{definition}[proposition]{Definition}
\newtheorem{notation}[proposition]{Notation}

\theoremstyle{remark}
\newtheorem{remark}[proposition]{Remark}

\numberwithin{equation}{section}

\begin{document}

\title{Existence in critical spaces for the magnetohydrodynamical system\\in 3D bounded Lipschitz domains}
\author{Sylvie Monniaux\footnote{Aix-Marseille Univ., CNRS, Centrale Marseille, I2M UMR7373, Marseille, France 
- {\tt sylvie.monniaux@univ-amu.fr}} \footnote{partially supported by the ANR project INFAMIE, ANR-15-CE40-0011}}
\date{}


\maketitle

\abstract{Existence of mild solutions for the 3D MHD system in bounded Lipschitz domains 
is established in critical spaces with the {\sl absolute} boundary conditions.}

\section{Introduction}

The magnetohydrodynamical system in a domain $\Omega\subset{\mathds{R}}^3$ on a 
time interval $(0,T)$ ($0<T\le \infty$) as considered in \cite{ST83} (with all constants equal to 1) 
reads
\begin{equation}
\label{mhd}\tag{MHD}
\left\{
\begin{array}{rclcl}
\partial_t u-\Delta u+\nabla \pi+(u\cdot\nabla)u&=&({\rm curl}\,b)\times b
&\mbox{ in }&(0,T)\times\Omega\\
\partial_t b-\Delta b&=&{\rm curl}\,(u\times b)&\mbox{ in }&(0,T)\times\Omega\\
{\rm div}\,u&=&0&\mbox{ in }&(0,T)\times\Omega\\
{\rm div}\,b&=&0&\mbox{ in }&(0,T)\times\Omega\\
\end{array}
\right.
\end{equation}
where $u:(0,T)\times\Omega\to{\mathds{R}}^3$ denotes the {\sl velocity} of the 
(incompressible homogeneous) fluid,
the {\sl magnetic field} (in the absence of magnetic monopole) is denoted by 
$b:(0,T)\times\Omega\to{\mathds{R}}^3$ and $\pi:(0,T)\times\Omega\to{\mathds{R}}^3$
is the pressure of the fluid. The first equation of \eqref{mhd} corresponds to {\sl Navier-Stokes}
equations subject to
the {\sl Laplace force} $({\rm curl}\,b)\times b$ applied by the magnetic field $b$. 
Actually, the divergence-free condition on the magnetic field $b$ comes from the fact that $b$ 
is in the range of the ${\rm curl}$ operator. The second equation of \eqref{mhd} 
describes the evolution of the magnetic field following the so-called {\sl induction} equation.

This system \eqref{mhd} (with $T=\infty$ and $\Omega={\mathds{R}}^3$) is invariant under 
the scaling 
$u_\lambda(t,x)=\lambda u(\lambda^2t,\lambda x)$, 
$b_\lambda(t,x)=\lambda b(\lambda^2t,\lambda x)$ and
$\pi_\lambda(t,x)=\lambda^2 \pi(\lambda^2t,\lambda x)$, $\lambda>0$. This suggests that 
a critical space for $(u,b)$
is ${\mathscr{C}}([0,\infty);L^3({\mathds{R}}^3)^3)\times 
{\mathscr{C}}([0,\infty);L^3({\mathds{R}}^3)^3)$.

The purpose of this paper is to prove existence of solutions of this system in this critical 
space in a bounded Lipschitz domain under the so-called 
{\sl absolute} boundary conditions, denoted by \eqref{bc1} below. This is investigated in 
Theorem~\ref{thm:mhd1global}, Theorem~\ref{thm:mhd1local} in Section~\ref{MHD1}. 
The methods used here come from the theory developed in \cite{McIM18} for the absolute 
boundary conditions.

In Section~\ref{tools} are collected results on potential operators (similar to the famous 
Bogovsk\u \i i operator),the Stokes operators with Dirichlet boundary conditions and Hodge 
boundary conditions, as well as properties of the Hodge Laplacian in bounded Lipschitz domains. 
Section~\ref{MHD1} is devoted to the existence of mild solutions of the system 
\eqref{mhd} under absolute boundary conditions on a bounded Lipschitz domain in critical 
spaces.

\section{Tools}
\label{tools}

In this section are recalled some results proved in \cite{McIM18} which will be useful in the 
following. See also \cite{MM09a} and \cite{MM09b}.

\begin{notation}
For an (unbounded) operator $A$ on a Banach space $X$, we denote by ${\rm{\sf D}}(A)$ 
its domain, ${\rm{\sf R}}(A)$ its range and ${\rm{\sf N}}(A)$ its null space.
\end{notation}

\subsection{Differential forms, Potential operators}

We consider the {\sl exterior derivative} $d:=\nabla\wedge=\sum_{j=1}^n \partial_j e_j\wedge$ and 
the {\sl interior deri\-va\-tive} (or co-derivative) 
$\delta:=-\nabla\lrcorner\,=-\sum_{j=1}^n \partial_j e_j\lrcorner\,$ acting on 
{\sl differential forms} on a domain $\Omega\subset{\mathbb{R}}^n$, i.e. 
acting on functions from $\Omega$ to  the exterior algebra 
$\Lambda=\Lambda^0\oplus\Lambda^1\oplus\dots\oplus\Lambda^n$ of
${\mathbb{R}}^n$. 

We denote by $\bigl\{e_S\,;\,S\subset\{1,\dots,n\}\bigr\}$ the basis
for $\Lambda$. The space of $\ell$-vectors $\Lambda^\ell$ is the
span of $\bigl\{e_S\,;\,|S|=\ell\bigr\}$, where 
\[
e_S=e_{j_1}\wedge e_{j_2}\wedge\dots\wedge e_{j_\ell}\quad
\mbox{for}\quad S=\{e_{j_1},\dots,e_{j_\ell}\} \quad\text{with }\ j_1<j_2<\dots<j_{\ell}.
\]
Remark that $\Lambda^0$, the space of complex scalars, is the span of 
$e_\emptyset$ ($\emptyset$ being the empty set). We set 
$\Lambda^\ell=\{0\}$ if $\ell<0$ or $\ell>n$.

On the exterior algebra $\Lambda$, the basic operations are
\begin{enumerate}[$(i)$ ]
\item 
the exterior product $\wedge:\Lambda^k\times\Lambda^\ell\to\Lambda^{k+\ell}$,
\item
the interior product $\lrcorner\,:\Lambda^k\times\Lambda^\ell\to\Lambda^{\ell-k}$,
\item
the Hodge star operator $\star:\Lambda^\ell\to\Lambda^{n-\ell}$, 
\item
the inner product $\langle\cdot,\cdot\rangle:\Lambda^\ell\times\Lambda^\ell
\to{\mathbb{R}}$. 
\end{enumerate}
If $a\in\Lambda^1$, $u\in\Lambda^\ell$ and $v\in\Lambda^{\ell+1}$, then
$$
\langle a\wedge u,v\rangle=\langle u,a\lrcorner\,v\rangle.
$$
For more details, we refer to, e.g., \cite[Section~2]{AMcI04} and 
\cite[Section~2]{CMcI10}, noting that both these papers contain some historical 
background (and being careful that $\delta$ has the opposite sign in 
\cite{AMcI04}). In particular, we note the relation between $d$ and
$\delta$ via the Hodge star operator:
\begin{equation}
\label{eq:*d=delta}
\star\delta u=(-1)^\ell d(\star\,u) \quad\mbox{and}\quad
\star du=(-1)^{\ell-1}\delta(\star\,u) \quad\mbox{for an $\ell$-form }u.  
\end{equation}
In dimension $n=3$, this gives (see \cite[\S2]{CMcI10}) for a vector $a\in{\mathds{R}}^3$ 
identified with a 1-form
\begin{itemize}
\item[-]
$u$ scalar, interpreted as 0-form: $a\wedge u=ua$, $a\lrcorner\, u=0$;
\item[-]
$u$ scalar, interpreted as 3-form: $a\wedge u=0$, $a\lrcorner\, u=ua$;
\item[-]
$u$ vector, interpreted as 1-form: $a\wedge u=a\times u$, $a\lrcorner\, u=a\cdot u$;
\item[-]
$u$ vector, interpreted as 2-form: $a\wedge u=a\cdot u$, $a\lrcorner\, u =-a\times u$.
\end{itemize}
The domains of the differential operators $d$ and $\delta$, denoted by 
${\rm{\sf D}}(d)$ and ${\rm{\sf D}}(\delta)$ are defined by
\[
{\rm{\sf D}}(d):=\bigl\{u\in L^2(\Omega,\Lambda); du\in L^2(\Omega,\Lambda)\bigr\}
\quad \mbox{and}\quad
{\rm{\sf D}}(\delta):=\bigl\{u\in L^2(\Omega,\Lambda); 
\delta u\in L^2(\Omega,\Lambda)\bigr\}.
\]
Similarly, the $L^p$ versions of these domains read
\[
{\rm{\sf D}}^p(d):=\bigl\{u\in L^p(\Omega,\Lambda); 
du\in L^p(\Omega,\Lambda)\bigr\}
\  \mbox{ and }\  
{\rm{\sf D}}^p(\delta):=\bigl\{u\in L^p(\Omega,\Lambda); 
\delta u\in L^p(\Omega,\Lambda)\bigr\}.
\]
The differential operators $d$ and $\delta$ satisfiy $d^2=d\circ d=0$ 
and $\delta^2=\delta\circ\delta=0$.
We will also consider the adjoints of $d$ and $\delta$ in the sense of
maximal adjoint operators in a Hilbert space: $\underline{\delta}:=d^*$
and $\underline{d}:=\delta^*$. They are defined as the closures in 
$L^2(\Omega,\Lambda)$ of the closable
operators $\bigl(d^*,{\mathscr{C}}_c^\infty(\Omega,\Lambda)\bigr)$ and
$\bigl(\delta^*,{\mathscr{C}}_c^\infty(\Omega,\Lambda)\bigr)$.

The following proposition has been proved in \cite[Proposition~4.1]{McIM18}
in a slightly more general framework (see also \cite[Theorem~1.5]{MMM08} and 
\cite[Theorem~1.1, Theorem~4.6 and Remark~4.12]{CMcI10}).

\begin{proposition}
\label{prop:propR,K}
Suppose $\Omega$ is a bounded Lipschitz  domain. 
Then the potential operators $R_\Omega$, $S_\Omega$ and $K_\Omega$ 
defined above satisfy for all $p\in (1,\infty)$, with the convention $p^S=\frac{np}{n-p}$ if $p<n$, 
$p^S=+\infty$ if $p>n$ and $p^S\in [n,+\infty)$ if $p=n$,
\begin{align*}
&R_\Omega: 
L^p(\Omega,\Lambda)\to L^{p^S}(\Omega,\Lambda)\cap{\textsf{D}}^p(d),
\quad 
S_\Omega:
L^p(\Omega,\Lambda)\to L^{p^S}(\Omega,\Lambda)\cap
{\textsf{D}}^p(d^*),
\\[4pt]
&K_\Omega:L^p(\Omega,\Lambda)\to L^\infty(\Omega,\Lambda)
\cap{\textsf{D}}^p(d), 
\quad
K_\Omega^*:L^p(\Omega,\Lambda)\to L^\infty(\Omega,\Lambda)
\cap{\textsf{D}}^p(d^*),
\\[4pt]
&K_\Omega, K_\Omega^*\mbox{ are compact operators in }L^p(\Omega,\Lambda),
\\[4pt]
&dR_\Omega+R_\Omega d={\rm I}\,-K_\Omega,
\qquad 
d^*S_\Omega+S_\Omega d^*={\rm I}\,-K_\Omega^*,
\\[4pt]
&dK_\Omega=0,\quad d^*K_\Omega^*=0 \quad \mbox{and}\quad
K_\Omega=0 \mbox{ on }{\rm{\sf R}}^p(d), \quad K_\Omega^*=0
\mbox{ on }{\rm{\sf R}}^p(d^*),
\\[4pt]
&dR_\Omega u=u\mbox{ if }u\in{\rm{\sf R}}^p(d), \quad
d^*S_\Omega u=u
\mbox{ if }u\in{\rm{\sf R}}^p(d^*).
\end{align*}
\end{proposition}

As direct consequence we obtain that $dR_\Omega$ and
$d^*S_\Omega$ are projections from 
$L^p(\Omega,\Lambda)$ onto the ranges of $d$ and $d^*$, ${\rm{\sf R}}^p(d)$ and 
${\rm{\sf R}}^p(d^*)$, for all $p\in(1,\infty)$. 

\subsection{Hodge-Laplacian and Hodge-Stokes operators in Lipschitz domains}

\begin{definition}
The {\sl Hodge-Dirac operator} on $\Omega$ with {\sl tangential boundary 
conditions} is
\[
D_{\|}:=d+d^*.
\]
Note that $-\Delta_{\|}:=D_{\|}^2=d d^*+ d^*d$ is 
the {\sl Hodge-Laplacian} with {\sl absolute} (generalised Neumann) boundary 
conditions.

For a scalar function 
$u:\Omega\to\Lambda^0$ we have that $-\Delta_{\|}u=d^*du=-\Delta_Nu$, where 
$\Delta_N$ is the Neumann Laplacian.
\end{definition}

Following \cite[Section 4]{AKMcI06Invent}, we have that the operator 
$D_{\|}$ is a closed densely defined operator
in $L^2(\Omega,\Lambda)$, and that 
\begin{align}
\label{H2}\tag{$H_2$}
L^2(\Omega,\Lambda)=&{\rm{\sf R}}(d)\stackrel{\bot}{\oplus}
{\rm{\sf R}}(d^*)\stackrel{\bot}{\oplus}{\rm{\sf N}}(D_{\|})\\
=&{\rm{\sf R}}(d)\stackrel{\bot}{\oplus}{\rm{\sf N}}(d^*)
\label{H2Rd}\\
=&{\rm{\sf N}}(d)\stackrel{\bot}{\oplus}{\rm{\sf R}}(d^*)
\label{H2Nd}
\end{align}
where ${\rm{\sf N}}(D_{\|})= {\rm{\sf N}}(d)\cap{\rm{\sf N}}(d^*)
={\rm{\sf N}}\bigl(\Delta_{\|}\bigr)$ is finite dimensional. The orthogonal 
projection from $L^2(\Omega,\Lambda)$ onto ${\rm{\sf N}}(d^*)$ (see \eqref{H2Rd}),
restricted to $1$-forms, is the 
well-known Helmholtz (or Leray) projection denoted by ${\mathbb{P}}$. Restricted
to $2$-forms, the orthogonal projection from $L^2(\Omega,\Lambda)$ onto 
${\rm{\sf R}}(d)$ will be denoted in the sequel by ${\mathbb{Q}}$. 

The $p$ version of the previous Hodge decompositions can be found in 
\cite[Theorem~4.3]{McIM18}: there exist Hodge 
exponents $p_H$, $p^H=p_H'$ with $1\le p_H<2<p^H\le \infty$ such that
\begin{align}
\label{Hp}\tag{$H_p$}
L^p(\Omega,\Lambda)=&{\rm{\sf R}}^p(d)\oplus
{\rm{\sf R}}^p(d^*)\oplus{\rm{\sf N}}(D_{\|})\\
=&{\rm{\sf R}}^p(d) \oplus{\rm{\sf N}}^p(d^*)\label{HpRd}\\
=&{\rm{\sf N}}^p(d) \oplus{\rm{\sf R}}^p(d^*)\label{HpNd}
\end{align}
for all $p\in(p_H,p^H)$ and the projections ${\mathbb{P}}:L^p(\Omega,\Lambda^1)\to 
{\rm{\sf N}}^p(d^*)_{|_{\Lambda^1}}$ and ${\mathbb{Q}}:L^p(\Omega,\Lambda^2)\to 
{\rm{\sf R}}^p(d)_{|_{\Lambda^2}}$ extend accordingly. 

\begin{remark}
\label{rk:pH}
If the domain is smooth or have a Lipschitz boundary, we have the following estimates
on the Hodge exponents $p_H$ and $p^H$.
\begin{enumerate}
\item
If $\Omega\subset {\mathds{R}}^n$ is smooth, then $p_H=1$ and $p^H=\infty$ 
(see \cite[Theorems~2.4.2 and 2.4.14]{Schw95}. 
\item
In the case of a bounded Lipschitz domain, 
$p_H<\frac{2n}{n+1}$ and consequently $p^H>\frac{2n}{n-1}$, which gives in dimension $n=3$: 
$p_H<\frac{3}{2}$ and $p^H>3$ (see \cite[\S7]{McIM18}).
\end{enumerate}
\end{remark}

\begin{remark}
\label{rk:PQ}
Proposition~\ref{prop:propR,K} and the projections ${\mathbb{P}}$ and ${\mathbb{Q}}$
yield
\begin{align*}
{\mathbb{P}}(R_\Omega du+K_\Omega u)=&\,u\quad\mbox{for }
u\in {\rm{\sf N}}^p(d^*)_{|_{\Lambda^1}},\\
{\mathbb{Q}}(S_\Omega d^*b+K^*_\Omega b)=&\,b\quad\mbox{for }
b\in{\rm{\sf R}}^p(d)_{|_{\Lambda^2}}
\end{align*}
for all $p_H<p<p^H$. The second equation comes from the fact that 
${\rm{\sf R}}^p(d)\subset {\rm{\sf N}}^p(d)$, using \eqref{HpNd}.
\end{remark}

The following results can be found partly in \cite[Theorem~7.3]{MM09a} (sectoriality)
and in \cite[\S8]{McIM18} (improvement of the interval of $p$ for the Hodge-Stokes operator
and bounded holomorphic functional calculus):

\begin{theorem}
\label{thm:HodgeL&S}
Suppose $\Omega$ is a bounded Lipschitz domain in ${\mathbb{R}}^n$.
Define $-\Delta_{\|}=D_{\|}^2$ in $L^2(\Omega,\Lambda)$. 
If $p_H<p<p^H$, then the Hodge-Laplacian with absolute boundary conditions
$-\Delta_{\|}$ is sectorial of angle $0$ in 
$L^p(\Omega,\Lambda)$ and for all $\mu\in (0,\frac{\pi}{2})$, $-\Delta_{\|}$
admits a bounded $S_{\mu +}^\circ$ holomorphic functional calculus in
$L^p(\Omega,\Lambda)$.

Define the Hodge-Stokes operator by $S_{\|}:=D_{\|}^2 = d^*d$ in 
${\rm{\sf N}}^2(d^*)$, restricted to $1$-forms. If $\max\bigl\{1,\frac{np_H}{n+p_H}\bigr\}<p<p^H$, 
then $S_{\|}$ is sectorial of angle $0$ in ${\rm{\sf N}}^p(d^*)_{|_{\Lambda^1}}$ 
and for all $\mu\in (0,\frac{\pi}{2})$, $S_{\|}$ admits a bounded 
$S_{\mu +}^\circ$ holomorphic functional calculus in 
${\rm{\sf N}}^p(d^*)_{|_{\Lambda^1}}$. In particular, the semigroup
$(e^{-tS_{\|}})_{t\ge 0}$ is bounded on ${\rm{\sf N}}^p(d^*)_{|_{\Lambda^1}}$ with 
norm denoted by $K_{p,S}$.

Define the Hodge-Maxwell operator $M_{\|}:=dd^*$ in ${\rm{\sf N}}(d)$, restricted to $2$-forms.
If $\max\bigl\{1,\frac{np_H}{n+p_H}\bigr\}<p<p^H$, then $M_{\|}$ is sectorial of angle $0$ in 
${\rm{\sf N}}^p(d)_{|_{\Lambda^2}}$ and for all $\mu\in (0,\frac{\pi}{2})$, $M_{\|}$ admits a 
bounded $S_{\mu +}^\circ$ holomorphic functional calculus in 
${\rm{\sf N}}^p(d)_{|_{\Lambda^2}}$.
In particular, the semigroup $(e^{-tM_{\|}})_{t\ge 0}$ is bounded on 
${\rm{\sf R}}^p(d)_{|_{\Lambda^1}}$ with norm denoted by $K_{p,M}$.
\end{theorem}

Using the results stated in Remark~\ref{rk:PQ}, one can prove $L^p-L^q$ bounds for the 
operator $S_{\|}$ (resp. $M_{\|}$) (see \cite[Theorems~3.1 and 4.1]{MM09b} for the dimension 3 
and \cite[Theorem~1.1]{HMM} for the Riesz transform like estimates \eqref{eq:RieszStokes} and
\eqref{eq:RieszMaxwell}).

\begin{theorem}
\label{thm:Lp-LqHodgeS}
Let $p\in\bigl(\max\bigl\{1,\frac{np_H}{n+p_H}\bigr\},p^H\bigr)$ and $q\in [p,p^H)$
such that $\frac{1}{p}-\frac{\alpha}{n}=\frac{1}{q}$ for some $\alpha\in[0,1]$. 
Then the semigroup $(e^{-tS_{\|}})_{t\ge 0}$ in ${\rm{\sf N}}^p(d^*)_{|_{\Lambda^1}}$
satisfies the estimates
\begin{equation}
\label{eq:Lp-LqStokes}
c^S_{p,q}:=\sup_{t\ge 0}\bigl\|t^{\frac{\alpha}{2}}
e^{-tS_{\|}}\bigr\|_{{\rm{\sf N}}^p(d^*)_{|_{\Lambda^1}}\to L^q}
+\sup_{t\ge 0}\bigl\|t^{\frac{1+\alpha}{2}}
d e^{-tS_{\|}}\bigr\|_{{\rm{\sf N}}^p(d^*)_{|_{\Lambda^1}}\to L^q}
<\infty
\end{equation}
and
\begin{equation}
\label{eq:RieszStokes}
\gamma_{p,q}^S:=\|S_{\|}^{-\frac{\alpha}{2}}\|_{{\rm{\sf N}}^p(d^*)_{|_{\Lambda^1}}\to L^q}<\infty.
\end{equation}
The semigroup $(e^{-tM_{\|}})_{t\ge 0}$ in ${\rm{\sf R}}^p(d)_{|_{\Lambda^2}}$
satisfies the estimate
\begin{equation}
\label{eq:Lp-LqMaxwell}
c^M_{p,q}:=\sup_{t\ge 0}\bigl\|t^{\frac{\alpha}{2}}
e^{-tM_{\|}}\bigr\|_{{\rm{\sf R}}^p(d)_{|_{\Lambda^2}}\to L^q}
+\sup_{t\ge 0}\bigl\|t^{\frac{1+\alpha}{2}}
d^* e^{-tM_{\|}}\bigr\|_{{\rm{\sf R}}^p(d)_{|_{\Lambda^2}}\to L^q}
<\infty
\end{equation}
and
\begin{equation}
\label{eq:RieszMaxwell}
\gamma_{p,q}^M:=\|M_{\|}^{-\frac{\alpha}{2}}\|_{{\rm{\sf R}}^p(d)_{|_{\Lambda^2}}\to L^q}<\infty.
\end{equation}
\end{theorem}

\section{Existence in the case of absolute boundary conditions}
\label{MHD1}

Thanks to the formula 
\[
(u\cdot\nabla)u=\tfrac{1}{2}\nabla|u|^2+u\times ({\rm curl}\,u)
\]
for a sufficiently smooth vector field $u$, the system \eqref{mhd} can be reformulated as follows:
\begin{equation}
\label{mhd1}
\left\{
\begin{array}{rclcl}
\partial_t u-\Delta u+\nabla \pi_1-u\times({\rm curl}\,u)&=&({\rm curl}\,b)\times b&\mbox{ in }&(0,T)\times\Omega\\
\partial_t b-\Delta b&=&{\rm curl}\,(u\times b)&\mbox{ in }&(0,T)\times\Omega\\
{\rm div}\,u&=&0&\mbox{ in }&(0,T)\times\Omega\\
{\rm div}\,b&=&0&\mbox{ in }&(0,T)\times\Omega\\
\end{array}
\right.
\end{equation}
where the pressure $\pi$ has been replaced by the so-called {\sl dynamical pressure} $\pi_1=\pi+\frac{1}{2}|u|^2$. 
This formulation can be translated in the language of differential forms: $\pi_1$ is a scalar function, interpreted as 
$0$-form, $u$ is a vector field interpreted as $1$-form and $b$ is a vector field interpreted as $2$-form. 
Following Section~\ref{tools} one can rewrite \eqref{mhd1} in terms of differential forms:
\begin{equation}
\label{mhd1diff}\tag{MHD1}
\left\{
\begin{array}{rclcl}
\partial_t u+S_{\|} u+d \pi_1+ u\lrcorner\,du&=&-d^*b\lrcorner\, b&\mbox{ in }&(0,T)\times\Omega\\
\partial_t b+M_{\|} b&=&-d(u\lrcorner\, b)&\mbox{ in }&(0,T)\times\Omega\\
u(t,\cdot)&\in&{\rm{\sf N}}(d^*)_{|_{\Lambda^1}}&\mbox{ for all }&t\in(0,T)\\
b(t,\cdot)&\in&{\rm{\sf R}}(d)_{|_{\Lambda^2}}&\mbox{ for all }&t\in(0,T).\\
\end{array}
\right.
\end{equation}
The terms in the first equation are all $1$-forms, in the second equation the terms are all $2$-forms.
The {\sl absolute} boundary conditions associated with the previous system \eqref{mhd1diff} are defined by the term $d^*$
in $-\Delta_{\|}=(dd^*+d^*d)$:
\begin{equation}
\label{bc1}\tag{BC1}
\left\{
\begin{array}{r}
\left.
\begin{array}{rclcl}
\nu\cdot u=\nu\lrcorner\, u &=&0&\mbox{ on }&(0,T)\times\partial\Omega\\
-\nu\times {\rm curl}\,u= \nu\lrcorner\, du&=&0&\mbox{ on }&(0,T)\times\partial\Omega\\
\end{array}
\right\} \mbox{ absolute b.c. for the 1-form }u\\
\left.\begin{array}{rclcl}
-\nu\times b=\nu\lrcorner\, b&=&0 &\mbox{ on }&(0,T)\times\partial\Omega\\
\nu\,{\rm div}\,b=\nu\lrcorner\, db&=&0 &\mbox{ on }&(0,T)\times\partial\Omega.\\
\end{array}
\right\} \mbox{ absolute b.c. for the 2-form }b
\end{array}
\right.
\end{equation}
This formulation can be used, for instance, to study the magnetohydrodynamical system in 
dimensions greater than or equal to 2 with the same theoretical tools. Let us point out that 
these boundary conditions are different to those usually investigated in magnetohydrodynamical
problems, starting with the paper \cite{ST83}; see also \cite{AB20}. The boundary conditions
\eqref{bc1} in the case of Navier-Stokes equations ({\it i.e.} for $b=0$) have been studied in 
\cite{MM09b}; see also \cite{M13} and \cite{MS18}.

\begin{remark}
The last condition in \eqref{bc1} is void since $b\in{\rm{\sf R}}(d)_{|_{\Lambda^2}}$: $db=0$ 
in all $\Omega$.
\end{remark}

\begin{definition}
\label{mildsol1}
Let $\Omega\subseteq {\mathds{R}}^3$.
A mild solution of the system \eqref{mhd1diff} with absolute boundary conditions \eqref{bc1} and 
initial conditions $u_0\in {\rm{\sf N}}(d^*)_{|_{\Lambda^1}}$ and $b_0 \in 
{\rm{\sf R}}(d)_{|_{\Lambda^2}}$ is a pair $(u,b)$ of vector fields satisfying
\begin{align}
\label{mildsolmhd1u}
u(t)=&e^{-tS_{\|}}u_0+\int_0^te^{-(t-s)S_{\|}}{\mathbb{P}}\bigl(-u(s)\lrcorner\, du(s)\bigr)\,{\rm d}s
+\int_0^te^{-(t-s)S_{\|}}{\mathbb{P}}\bigl(-d^*b(s)\lrcorner\,b(s)\bigr)\,{\rm d}s,
\\
\label{mildsolmhd1b}
b(t)=&e^{-t M_{\|}}b_0+\int_0^te^{-(t-s) M_{\|}}\Bigl(-d\bigl(u(s)\lrcorner\,b(s)\bigr)\Bigr)\,{\rm d}s.
\end{align}
\end{definition}

From now on, we assume the following technical
(Leibniz rule-like) property on the domain $\Omega\subset{\mathds{R}}^3$: for all $q\in[3,p^H)$,
there exists a constant $C_q>0$ such that 
\begin{equation}
\label{Leibniz}
\|d(\omega_1\lrcorner\,\omega_2)\|_{\frac{q}{2}}\le 
C_q \bigl(\|D_{\|}\omega_1\|_q\|\omega_2\|_q +\|\omega_1\|_q\|D_{\|}\omega_2\|_q\bigr)
\end{equation}
for all $\omega_1\in {\rm{\sf D}}^p(D_{\|})\cap L^q(\Omega,\Lambda^1)$ and all
$\omega_2\in {\rm{\sf D}}^p(D_{\|})\cap L^q(\Omega,\Lambda^2)$. This is the case if
the domain $\Omega$ is smooth.

The following theorem is about the global existence of mild solutions with small initial data.

\begin{theorem}[Global existence]
\label{thm:mhd1global}
Let $\Omega\subset{\mathds{R}}^3$ be a bounded Lipschitz domain or 
$\Omega={\mathds{R}}^3$. Then there
exists $\varepsilon>0$ such that for all $u_0\in {\rm{\sf N}}^3(d^*)_{|_{\Lambda^1}}$ and 
$b_0\in {\rm{\sf R}}^3(d)_{|_{\Lambda^2}}$ with $\|u_0\|_3+\|b_0\|_3\le \varepsilon$, the
system \eqref{mhd1diff} with the boundary conditions \eqref{bc1} and $T=\infty$ admits a 
mild solution 
$u,b\in{\mathscr{C}}([0,\infty);L^3(\Omega)^3)$.
\end{theorem}

The next result states local existence of mild solutions with no restriction on the size of the 
initial data.

\begin{theorem}[Local existence]
\label{thm:mhd1local}
Let $\Omega\subset{\mathds{R}}^3$ be a bounded Lipschitz domain or 
$\Omega={\mathds{R}}^3$. Then for all 
$u_0\in {\rm{\sf N}}^3(d^*)_{|_{\Lambda^1}}$ and 
$b_0\in {\rm{\sf R}}^3(d)_{|_{\Lambda^2}}$ there exists $T>0$ such that the
system \eqref{mhd1diff} with the boundary conditions \eqref{bc1} admits a mild solution 
$u,b\in{\mathscr{C}}([0,T);L^3(\Omega)^3)$.
\end{theorem}

The methods to prove these two theorems are classical based on a fixed point theorem, 
already used for the Navier-Stokes equations in the paper by Fujita and Kato \cite{FK64}
(see also \cite{M06}) and in \cite{BM20} (see also \cite{BH20}) for the Boussinesq system. 
Most of the tools used here appeared in the paper \cite{MM09b}; see also \cite{McIM18}.

Let $q\in\bigl(3,\min\{p^H,6\}\bigr)$ and $\alpha\in(0,1)$ such that 
$\frac{1}{q}=\frac{1}{3}-\frac{\alpha}{3}$.
For $0<T\le\infty$, we define the following spaces
\begin{align}
\label{eq:UT}
{\mathscr{U}}_T:=&\bigl\{u\in{\mathscr{C}}((0,T);{\rm{\sf N}}^q(d^*)_{|_{\Lambda^1}});
du\in{\mathscr{C}}((0,T);L^q(\Omega,\Lambda^2)) :\\
\nonumber
&\sup_{0<t<T}\bigl(t^{\frac{\alpha}{2}}\|u(t)\|_q+t^{\frac{1+\alpha}{2}}\|du(t)\|_q\bigr)<\infty\bigr\}
\end{align}
and
\begin{align}
\label{eq:BT}
{\mathscr{B}}_T:=&\bigl\{b\in{\mathscr{C}}((0,T);{\rm{\sf R}}^q(d)_{|_{\Lambda^2}}) ;
d^*b\in{\mathscr{C}}((0,T);L^q(\Omega,\Lambda^1)) :\\
\nonumber
&\sup_{0<t<T}\bigl(t^{\frac{\alpha}{2}}\|b(t)\|_q+t^{\frac{1+\alpha}{2}}\|d^*b(t)\|_q\bigr)<\infty\bigr\},
\end{align}
endowed with the norms
\begin{equation}
\label{eq:normUT}
\|u\|_{{\mathscr{U}}_T}:=\sup_{0<t<T}\bigl(t^{\frac{\alpha}{2}}\|u(t)\|_q
+t^{\frac{1+\alpha}{2}}\|du(t)\|_q\bigr)
\end{equation}
and 
\begin{equation}
\label{eq:normBT}
\|b\|_{{\mathscr{B}}_T}:=\sup_{0<t<T}\bigl(t^{\frac{\alpha}{2}}\|b(t)\|_q
+t^{\frac{1+\alpha}{2}}\|d^*b(t)\|_q\bigr).
\end{equation}

\begin{lemma}
\label{lem:initcond}
For $u_0\in {\rm{\sf N}}^3(d^*)_{|_{\Lambda^1}}$ and 
$b_0\in {\rm{\sf R}}^3(d)_{|_{\Lambda^2}}$, we have
\begin{enumerate}
\item
$a_1:t\mapsto e^{-tS_{\|}}u_0\in{\mathscr{U}}_T$,
\item
$a_2:t\mapsto e^{-tM_{\|}}b_0\in{\mathscr{B}}_T$,
\end{enumerate}
for all $T>0$
Moreover, for all $\varepsilon>0$, there exists $T>0$ such that
\begin{equation}
\label{eq:a1a2}
\|a_1\|_{{\mathscr{U}}_T}+\|a_2\|_{{\mathscr{B}}_T}\le \varepsilon.
\end{equation}
\end{lemma}

\begin{proof}
By Theorem~\ref{thm:Lp-LqHodgeS}, the following bound holds for all $T>0$:
 \begin{equation}
 \label{eq:esta1a2}
 \|a_1\|_{{\mathscr{U}}_T}+\|a_2\|_{{\mathscr{B}}_T}\le 
c_{3,q}^S \|u_0\|_3 + c^M_{3,q}\|b_0\|_3.
 \end{equation}
 Therefore, if $\|u_0\|_3$ and $\|b_0\|_3$ are small enough, \eqref{eq:a1a2} holds for
 every $T>0$.
 
 For any $u_0$ and $b_0$ (not necessarily small in the $L^3$ norm), for $\epsilon>0$,
 let $u_{0,\epsilon}\in {\rm{\sf N}}^q(d^*)_{|_{\Lambda^1}}$ and 
 $b_{0,\epsilon}\in {\rm{\sf R}}^3(d)_{|_{\Lambda^2}}$ such that
 \[
 \|u_0-u_{0,\epsilon}\|_3+\|b_0-b_{0,\epsilon}\|_3\le \epsilon.
 \]
 We denote by $a_{1,\epsilon}$ and $a_{2,\epsilon}$ the quantities
 $a_{1,\epsilon}(t)=e^{-tS_{\|}}u_{0,\epsilon}$ and 
 $a_{2,\epsilon}(t)=e^{-tM_{\|}}b_{0,\epsilon}$. By \eqref{eq:esta1a2},
 there holds
 \begin{equation}
 \label{eq:a1epsa2eps}
 \|a_1-a_{1,\epsilon}\|_{{\mathscr{U}}_T}+\|a_2-a_{2,\epsilon}\|_{{\mathscr{B}}_T}\le 
\epsilon\,(c_{3,q}^S+ c^M_{3,q}).
 \end{equation}
 Applying \eqref{eq:Lp-LqStokes} with $p=q$ and $\alpha=0$, we obtain
 \[
 \|a_{1,\epsilon}\|_{{\mathscr{U}}_T}\le c^s_{q,q} T^{\frac{\alpha}{2}}\|u_{0,\epsilon}\|_q.
 \]
 The same reasoning applying \eqref{eq:Lp-LqMaxwell} with $p=q$ and $\alpha=0$ yields
 \[
 \|a_{2,\epsilon}\|_{{\mathscr{B}}_T} \le c^M_{q,q} T^{\frac{\alpha}{2}}\|b_{0,\epsilon}\|_q.
 \]
 Now choosing $\epsilon>0$ small enough and $T>0$ small enough, 
 we find that \eqref{eq:a1a2} holds.
\end{proof}

Next, we define the operators
\begin{equation}
\label{eq:B1}
B_1(u,v)(t)=\int_0^te^{-(t-s)S_{\|}}{\mathbb{P}}\bigl(-u(s)\lrcorner\, dv(s)\bigr)\,{\rm d}s,
\quad t\in[0,T),\ u,v \in {\mathscr{U}}_T,
\end{equation}
\begin{equation}
\label{eq:B2}
B_2(b,b')(t)=\int_0^te^{-(t-s)S_{\|}}{\mathbb{P}}\bigl(-d^*b(s)\lrcorner\, b'(s)\bigr)\,{\rm d}s,
\quad t\in[0,T),\ b,b'\in {\mathscr{B}}_T, 
\end{equation}
\begin{equation}
\label{eq:B3}
B_3(u,b)(t)=\int_0^te^{-(t-s)M_{\|}}\Bigl(-d\bigl(u(s)\lrcorner\, b(s)\bigr)\Bigr)\,{\rm d}s,
\quad t\in[0,T),\ u \in {\mathscr{U}}_T,\ b\in {\mathscr{B}}_T. 
\end{equation}
The next lemma gives a precise statement about the boundedness of the bilinear operators 
$B_1$, $B_2$ and $B_3$.

\begin{lemma}
\label{lem:B1B2B3}
The bilinear operators $B_1$, $B_2$ and $B_3$ are bounded in the following spaces:
\begin{enumerate}
\item
$B_1:{\mathscr{U}}_T\times{\mathscr{U}}_T\to{\mathscr{U}}_T$,
\item
$B_2:{\mathscr{B}}_T\times{\mathscr{B}}_T\to{\mathscr{U}}_T$,
\item
$B_3:{\mathscr{U}}_T\times{\mathscr{B}}_T\to{\mathscr{B}}_T$
\end{enumerate}
with norms independent from $T>0$.
\end{lemma}

\begin{proof}
\begin{enumerate}
\item
For $u,v\in{\mathscr{U}}_T$, by definition of ${\mathscr{U}}_T$ we have that
$s\mapsto s^{\frac{1}{2}+\alpha}u(s)\lrcorner\,dv(s) 
\in {\mathscr{C}}_b((0,T);L^{\frac{q}{2}}(\Omega,\Lambda^1)$
with norm less than or equal to $\|u\|_{{\mathscr{U}}_T}\|v\|_{{\mathscr{U}}_T}$.
Since $\frac{3}{2}<\tfrac{q}{2}<3$, ${\mathbb{P}}$ is bounded from 
$L^{\frac{q}{2}}(\Omega,\Lambda^1)$ to ${\rm{\sf N}}^{\frac{q}{2}}(d^*)_{|_{\Lambda^1}}$
Moreover, $e^{-(t-s)S_{\|}}$ maps ${\rm{\sf N}}^{\frac{q}{2}}(d^*)_{|_{\Lambda^1}}$
to ${\rm{\sf N}}^q(d^*)_{|_{\Lambda^1}}$ with norm $c_{\frac{q}{2},q}^S (t-s)^{-\frac{1-\alpha}{2}}$
thanks to \eqref{eq:Lp-LqStokes} with $p=\frac{q}{2}$. Therefore, we have for all $t\in(0,T)$
\[
\|B_1(u,v)(t)\|_q\lesssim 
\Bigl(\int_0^t s^{-\frac{1}{2}-\alpha}(t-s)^{-\frac{1-\alpha}{2}}\,{\rm d}s\Bigr)
\|u\|_{{\mathscr{U}}_T}\|v\|_{{\mathscr{U}}_T}
\lesssim t^{-\frac{\alpha}{2}}\|u\|_{{\mathscr{U}}_T}\|v\|_{{\mathscr{U}}_T}.
\] 
This gives the first estimate for $B_1(u,v)\in{\mathscr{U}}_T$. For the second estimate,
we note that $de^{-(t-s)S_{\|}}$ maps ${\rm{\sf N}}^{\frac{q}{2}}(d^*)_{|_{\Lambda^1}}$
to $L^q(\Omega,\Lambda^2)$ with norm $c_{\frac{q}{2},q}^S (t-s)^{-1+\frac{\alpha}{2}}$
thanks to \eqref{eq:Lp-LqStokes} with $p=\frac{q}{2}$. Therefore, we have for all $t\in(0,T)$
\[
\|dB_1(u,v)(t)\|_q\lesssim 
\Bigl(\int_0^t s^{-\frac{1}{2}-\alpha}(t-s)^{-1+\frac{\alpha}{2}}\,{\rm d}s\Bigr)
\|u\|_{{\mathscr{U}}_T}\|v\|_{{\mathscr{U}}_T}
\lesssim t^{-\frac{1+\alpha}{2}}\|u\|_{{\mathscr{U}}_T}\|v\|_{{\mathscr{U}}_T},
\] 
which gives the second estimate for $B_1(u,v)\in{\mathscr{U}}_T$.
\item
The proof that for $b,b'\in{\mathscr{B}}_T$, $B_2(b,b')\in {\mathscr{U}}_T$ with norm 
independent from $T>0$ follows the lines of the previous point. We omit the details here.
\item
Thanks to the property \eqref{Leibniz}, the proof that for $u\in{\mathscr{U}}_T$ and
$b\in{\mathscr{B}}_T$, $B_3(u,b)\in{\mathscr{B}}_T$ with norm independent 
from $T>0$ can be copied from the proof of point 1, using the fact that 
$d^*e^{-(t-s)M_{\|}}$ maps ${\rm{\sf R}}^{\frac{q}{2}}(d)_{|_{\Lambda^2}}$ to
$L^q(\Omega,\Lambda^1)$ with norm $c_{\frac{q}{2},q}^M(t-s)^{-1+\frac{\alpha}{2}}$
thanks to \eqref{eq:Lp-LqMaxwell} with $p=\frac{q}{2}$.
\end{enumerate}
This proves Lemma~\ref{lem:B1B2B3}.
\end{proof}

\begin{lemma}
\label{lem:UTBTmild}
Let $T>0$.
Assume that $(u,b)\in {\mathscr{U}}_T\times{\mathscr{B}}_T$ is a mild solution of
\eqref{mhd1diff} with absolute boundary conditions \eqref{bc1} with initial conditions 
$u_0\in {\rm{\sf N}}^3(d^*)_{|_{\Lambda^1}}$ and 
$b_0\in {\rm{\sf R}}^3(d)_{|_{\Lambda^2}}$. Then
$u\in {\mathscr{C}}_b([0,T);{\rm{\sf N}}^3(d^*)_{|_{\Lambda^1}})$ and
$b\in {\mathscr{C}}_b([0,T);{\rm{\sf R}}^3(d)_{|_{\Lambda^2}})$.
\end{lemma}

\begin{proof}
To prove this lemma, first observe that if $u_0\in {\rm{\sf N}}^3(d^*)_{|_{\Lambda^1}}$ and 
$b_0\in {\rm{\sf R}}^3(d)_{|_{\Lambda^2}}$, then for all $T>0$,
$t\mapsto e^{-tS_{\|}}u_0 \in {\mathscr{C}}_b([0,T);{\rm{\sf N}}^3(d^*)_{|_{\Lambda^1}})$
and $t\mapsto e^{-tM_{\|}}b_0\in  {\mathscr{C}}_b([0,T);{\rm{\sf R}}^3(d)_{|_{\Lambda^2}})$.
It remains to show that if $u\in{\mathscr{U}}_T$ and $b\in{\mathscr{B}}_T$, then
$B_1(u,u)\in {\mathscr{C}}_b([0,T);{\rm{\sf N}}^3(d^*)_{|_{\Lambda^1}})$, 
$B_2(b,b)\in {\mathscr{C}}_b([0,T);{\rm{\sf N}}^3(d^*)_{|_{\Lambda^1}})$
and $B_3(u,b)\in {\mathscr{C}}_b([0,T);{\rm{\sf R}}^3(d)_{|_{\Lambda^2}})$.
The continuity is straightforward. To prove boundedness, it suffices to reproduce the proof of 
the previous lemma (recall that $\alpha=1-\frac{3}{q}$) to obtain
\[
\|B_1(u,u)(t)\|_3\lesssim 
\Bigl(\int_0^t s^{-\frac{1}{2}-\alpha}(t-s)^{-\frac{1}{2}+\alpha}\,{\rm d}s\Bigr)
\|u\|_{{\mathscr{U}}_T}^2
\lesssim \|u\|_{{\mathscr{U}}_T}^2,
\]
using the fact that $e^{-(t-s)S_{\|}}$ maps ${\rm{\sf N}}^{\frac{q}{2}}(d^*)_{|_{\Lambda^1}}$ 
to $L^3(\Omega,\Lambda^1)$ with norm controlled by 
$c_{\frac{q}{2},3}^S(t-s)^{-\frac{3}{q}+\frac{1}{2}}$ thanks to \eqref{eq:Lp-LqStokes}.
The terms $B_2$ and $B_3$ can be treated similarly.
\end{proof}

\begin{proof}[Proof of Theorems~\ref{thm:mhd1global} and \ref{thm:mhd1local}]
The system
\begin{align}
\label{eq:fixedpoint}
u=a_1+B_1(u,u)+B_2(b,b)\quad\mbox{and}\quad b=a_2+B_3(u,b), \quad (u,b)\in{\mathscr{U}}_T
\end{align}
can be reformulated as
\begin{equation}
\label{eq:picard}
{\boldsymbol{u}}={\boldsymbol{a}}+{\boldsymbol{B}}({\boldsymbol{u}},{\boldsymbol{u}})
\end{equation}
where ${\boldsymbol{u}}=(u,b)\in {\mathscr{U}}_T\times{\mathscr{B}}_T$, ${\boldsymbol{a}}=(a_1,a_2)$ and 
${\boldsymbol{B}}({\boldsymbol{u}},{\boldsymbol{v}})=(B_1(u,v)+B_2(b,b'),B_3(u,b'))$ if 
${\boldsymbol{u}}=(u,b)$ and ${\boldsymbol{v}}=(v,b')$. On ${\mathscr{U}}_T\times{\mathscr{B}}_T$
we choose the norm 
$\|(u,b)\|_{{\mathscr{U}}_T\times{\mathscr{B}}_T}:=\|u\|_{{\mathscr{U}}_T}+\|b\|_{{\mathscr{B}}_T}$.
One can easily check, using Lemma~\ref{lem:B1B2B3}, that
\[
\|{\boldsymbol{B}}({\boldsymbol{u}},{\boldsymbol{v}})\|_{{\mathscr{U}}_T\times{\mathscr{B}}_T}
\le C \|{\boldsymbol{u}}\|_{{\mathscr{U}}_T\times{\mathscr{B}}_T} 
\|{\boldsymbol{v}}\|_{{\mathscr{U}}_T\times{\mathscr{B}}_T}
\]
where $C$ is a constant independent from $T>0$. We can then apply Picard's fixed point theorem
to prove that for $u_0\in {\rm{\sf N}}^3(d^*)_{|_{\Lambda^1}}$ and 
$b_0\in {\rm{\sf R}}^3(d)_{|_{\Lambda^2}}$, with $T\le \infty$ such that 
\eqref{eq:a1a2} holds for $\varepsilon=\frac{1}{4C}$, the system \eqref{eq:picard} 
admits a unique solution ${\boldsymbol{u}}=(u,b)\in {\mathscr{U}}_T\times{\mathscr{B}}_T$.
By Lemma~\ref{lem:UTBTmild}, this provides a mild solution 
$(u,b)\in{\mathscr{C}}_b([0,T);{\rm{\sf N}}^3(d^*)_{|_{\Lambda^1}})\times
{\mathscr{C}}_b([0,T);{\rm{\sf R}}^3(d)_{|_{\Lambda^2}})$
of \eqref{mhd1diff} with boundary conditions \eqref{bc1}.
\end{proof}



\end{document}